\input amstex
\documentstyle{amsppt}
%
%
\nopagenumbers
\def\const{\operatorname{const}}
\def\negskp{\hskip -2pt}
\def\compos{\,\raise 1pt\hbox{$\sssize\circ$} \,}
\pagewidth{360pt}
\pageheight{606pt}
\rightheadtext{Dynamical systems admitting normal shift and \dots}
\topmatter
\title Dynamical systems admitting normal shift
and wave equations.
\endtitle
\author
R.~A.~Sharipov
\endauthor
\abstract
High frequency limit for most of wave phenomena is known
as quasiclassical limit or ray optics limit. Propagation
of waves in this limit is described in terms of wave fronts
and rays. Wave front is a surface of constant phase whose
points are moving along rays. As it appears, their motion
can be described by Hamilton equations being special case
for Newton's equations. In simplest cases (e\.\,g\. light
propagation in isotropic non-homogeneous refracting media)
the dynamics of points preserves orthogonality of wave
front and rays. This property was generalized in the theory
of dynamical systems admitting normal shift of hypersurfaces.
In present paper inverse problem is studied, i\.\,e\. the
problem of reconstructing wave equation corresponding wave
front dynamics for which is described by Newtonian dynamical
system admitting normal shift of hypersurfaces in Riemannian
manifold.
\endabstract
\address Rabochaya street 5, 450003, Ufa, Russia
\endaddress
\email \vtop to 20pt{\hsize=280pt\noindent
R\_\hskip 1pt Sharipov\@ic.bashedu.ru\newline
r-sharipov\@mail.ru\vss}
\endemail
\urladdr
http:/\negskp/www.geocities.com/CapeCanaveral/Lab/5341
\endurladdr
\endtopmatter
\loadbold
\TagsOnRight
\document
\head
1. High frequency asymptotics for wave equations.
\endhead
    Let $M$ be a Riemannian manifold of the dimension $n$.
Denote by $\nabla$ covariant differentiation with respect
to metric connection in $M$. Then the expression
$$
\hskip -2em
\bold a_r\cdot D^r=\sum^n_{k_1=1}\!...\!\sum^n_{k_r\,=1}
\left(\frac{1}{i}\right)^{\!r}
a^{k_1\ldots\,k_r}(x^1,\ldots,x^n)\cdot
\nabla_{k_1}\cdot\ldots\cdot\nabla_{k_r}
\tag1.1
$$
is an elementary differential operator of $r$-th order. Here 
$a^{k_1\ldots\,k_r}(x^1,\ldots,x^n)$ are components of some
smooth symmetric tensor field $\bold a$ of the type $(r,0)$
in local coordinates $x^1,\,\ldots,\,x^n$ and $i=\sqrt{-1}$.
Operators \thetag{1.1} are united into polynomial
$$
\hskip -2em
H(p,D)=\sum^m_{r=0}\bold a_r\cdot D^r.
\tag1.2
$$
This is $m$-th order scalar differential operator in $M$. Here
$p$ is a point of $M$ and $D$ is a formal symbol for differentiation.
Operator \thetag{1.2} can be applied either to scalar field or
tensorial field in $M$. We shall apply it to scalar field $u$,
but first we introduce large parameter $\lambda$ in $H(p,D)$. Let's
denote
$$
\hskip -2em
H(p,\lambda^{-1}D)=\sum^m_{r=0}\lambda^{-r}\,\bold a_r\cdot D^r.
\tag1.3
$$
Then consider scalar $m$-th order wave equation with large parameter
$\lambda$:
$$
\hskip -2em
H(p,\lambda^{-1}D)u=0.
\tag1.4
$$
In particular cases \thetag{1.4} can reduce to Helmholtz equation,
to Schredinger equation, and to wave equation (see \cite{1}). The
following ansatz for formal asymptotic solution of the equation
\thetag{1.4} as $\lambda\to\infty$ was suggested by P.~Debye:
$$
\hskip -2em
u=\sum^\infty_{\alpha=0}\frac{\varphi_{\sssize(\ssize\alpha\sssize)}}
{(i\,\lambda)^\alpha}\cdot e^{i\lambda S}.
\tag1.5
$$
Here $\varphi_{\sssize(\ssize\alpha\sssize)}$ and $S$ are some
smooth scalar fields in $M$. Scalar field $S$ is interpreted as
a phase of wave.\par
    Now let's substitute \thetag{1.5} into \thetag{1.4}. This is
done in several steps. First of all let's apply differential
operator $\nabla_{k_1}\ldots\nabla_{k_r}$ to $e^{i\lambda S}$.
This yields
$$
\hskip -2em
\nabla_{k_1}\ldots\nabla_{k_r}e^{i\lambda S}=
\left(\sum^r_{q=0}\beta^{\sssize (\ssize q\sssize)}_{k_1\ldots\,k_r}
\cdot(i\,\lambda)^{r-q}\right)\cdot e^{i\lambda S}.
\tag1.6
$$
Leading coefficient in this expression can be written explicitly
$$
\hskip -2em
\beta^{\sssize (\lower 1pt\hbox{$\ssize 0$}\sssize)}_{k_1\ldots\,
k_r}=\nabla_{k_1}S\cdot\ldots\cdot\nabla_{k_r}S.
\tag1.7
$$
Other coefficients $\beta^{\sssize (\ssize q\sssize)}_{k_1\ldots\,
k_r}$ can also be calculated for each particular $q$. It might be
difficult to write general formula for all of them, but we do not
need this formula.\par
    In the second step we apply operator \thetag{1.3} to power series
\thetag{1.5}. Using formula \thetag{1.6}, for $H(p,\lambda^{-1}D)u$
we obtain the following expression:
$$
Hu=\sum^\infty_{\alpha=0}\sum^m_{r=0}
\sum^r_{s=0}\sum^{r}_{q=s}\sum^n_{k_1=1}\!...\!
\sum^n_{k_r\,=1}\frac{r!\cdot a^{k_1\ldots\,k_r}\cdot\beta^{\sssize
(\ssize q-s\kern -0.2pt\sssize)}_{k_{s+1}\ldots\,k_r}}{s!\,(r-s)!\cdot
(i\,\lambda)^{q+\alpha}}\cdot e^{i\lambda S}\cdot 
\nabla_{k_1}\!\ldots\!\nabla_{k_s}\varphi_{\sssize(\ssize\alpha
\sssize)}.
$$
If we collect all terms with fixed value of $\alpha$ and $\theta=q+\alpha$,
we find that the number of such terms in the above expression is finite.
Simplest way to prove this fact is to make the following transformations
of sums:
$$
\gather
\sum^\infty_{\alpha=0}\sum^m_{r=0}\sum^r_{s=0}\sum^{r}_{q=s}\to
\sum^\infty_{\alpha=0}\sum^m_{r=0}\sum^r_{q=0}\sum^{q}_{s=0}\to
\sum^\infty_{\alpha=0}\sum^m_{q=0}\sum^m_{r=q}\sum^{q}_{s=0}\ ,\\
\sum^\infty_{\alpha=0}\sum^m_{q=0}\sum^m_{r=q}\sum^{q}_{s=0}\to
\sum^\infty_{\theta=0}\ 
\sum^\theta_{\alpha=\alpha_{\lower 1pt\hbox{$\sssize\theta\!,\!m$}}}
\ \sum^m_{r=\theta-\alpha}\ \sum^{\theta-\alpha}_{s=0}\ .
\endgather
$$
Here $\alpha_{\theta,m}=\max(0,\theta-m)$. Now we introduce
differential operator $R_\omega$:
$$
\hskip -2em
R_\omega=\sum^m_{r=\omega}\ \sum^{\omega}_{s=0}\ 
\sum^n_{k_1=1}\!...\!\sum^n_{k_r\,=1}\frac{r!\cdot a^{k_1\ldots\,k_r}
\cdot\beta^{\sssize(\kern 0.3pt\ssize\omega-s\kern -0.2pt
\sssize)}_{k_{s+1}\ldots\,k_r}}{s!\,(r-s)!}\cdot
\nabla_{k_1}\!\ldots\!\nabla_{k_s}.
\tag1.8
$$
Let's rewrite left hand side of the equation \thetag{1.4}
in terms of operator \thetag{1.8}:
$$
Hu=\sum^\infty_{\theta=0}\ \sum^\theta_{\alpha=\alpha_{\lower 1pt
\hbox{$\sssize\theta\!,\!m$}}}\frac{R_{\theta-\alpha}\varphi_{\sssize(
\ssize\alpha\sssize)}}{(i\,\lambda)^\theta}\cdot e^{i\lambda S}.
$$
This is third step in substituting formal power series \thetag{1.5}
into the equation \thetag{1.4}. As a result wave equation \thetag{1.4}
breaks into infinite series of equations for $\varphi_{\sssize(\ssize
\alpha\sssize)}$:
$$
\hskip -2em
\sum^\theta_{\alpha=\alpha_{\lower 1pt\hbox{$\sssize\theta\!,\!m$}}}
R_{\theta-\alpha}\varphi_{\sssize(\ssize\alpha\sssize)}=0\text{\
\ for \ }\theta=0,\,1,\,\ldots,\,\infty.
\tag1.9
$$
In order to study the equations \thetag{1.9} we calculate explicit
form of operator \thetag{1.8} for $\omega=0$. Using
formula \thetag{1.7}, we derive
$$
\hskip -2em
R_0=\sum^m_{r=0}\sum^n_{k_1=1}\!...\!\sum^n_{k_r\,=1}
a^{k_1\ldots\,k_r}(x^1,\ldots,x^n)\cdot\nabla_{k_1}S
\cdot\ldots\cdot\nabla_{k_r}S.
\tag1.10
$$
Note that $R_0$ is a scalar operator, i\.\,e\. it is the operator
of multiplication by the function $R_0=R_0(x^1,\ldots,x^n)$.
Therefore first equation in the series \thetag{1.9} looks like
$R_0\cdot\varphi_{\sssize(\ssize\lower 1pt\hbox{$\ssize 0$}
\sssize)}=0$. Leading term in power series \thetag{1.5} for $u$
is nonzero. Therefore the equation $R_0\cdot\varphi_{\sssize(
\ssize\lower 1pt\hbox{$\ssize 0$}\sssize)}=0$ reduces to the
following one:
$$
\hskip -2em
R_0=0.
\tag1.11
$$
Taking into account \thetag{1.11}, we can rewrite \thetag{1.9} as
$$
\hskip -2em
\sum^{\theta-1}_{\alpha=\alpha_{\lower 1pt\hbox{$\sssize\theta\!,\!m$}}}
R_{\theta-\alpha}\varphi_{\sssize(\ssize\alpha\sssize)}=0\text{,
\ where \ }\theta=1,\,2,\,\ldots,\,\infty.
\tag1.12
$$
First equation in the series \thetag{1.12} is also very simple.
It is a partial differential equation of the first order for
the function $\varphi_{\sssize(\ssize\lower 1pt\hbox{$\ssize 0$}
\sssize)}$:
$$
\hskip -2em
R_1\varphi_{\sssize(\ssize\lower 1pt\hbox{$\ssize 0$}
\sssize)}=0.
\tag1.13
$$
Other equations in the series \thetag{1.12} can be written in
recursive form:
$$
\hskip -2em
R_1\varphi_{\sssize(\ssize\lower 1pt\hbox{$\ssize\theta-1$}
\sssize)}=-\sum^{\theta-2}_{\alpha=\alpha_{\lower 1pt\hbox{$\sssize
\theta\!,\!m$}}}R_{\theta-\alpha}\varphi_{\sssize(\ssize\alpha
\sssize)}\text{\ \ for \ }\theta=2,\,3,\,\ldots,\,\infty.
\tag1.14
$$
Below we calculate explicit expression for operator $R_1$, but we do
not consider the problem of solvability of equations \thetag{1.13}
and \thetag{1.14}. This should be done for each particular equation
\thetag{1.4} taking into account initial-value and/or boundary-value
problem stated for this equation.\par
    In order to find explicit expression for operator $R_1$ we should
first find explicit expression for coefficient $\beta^{\sssize (
\lower 1pt\hbox{$\ssize 1$}\sssize)}_{k_1\ldots\,k_r}$ in 
\thetag{1.6}. \pagebreak This is done by direct calculations. Applying
operator $\nabla_{k_1}\ldots\nabla_{k_r}$ to $e^{i\lambda S}$ and
collecting all terms with $(i\,\lambda)^{r-1}$, we get
$$
\gather
\hskip -2em
\beta^{\sssize(\lower 1pt\hbox{$\ssize 1$}\sssize)}_{k_1}
=0\text{\ \ \ for \ \ }r=1,
\tag1.15\\
\vspace{1ex}
\hskip -2em
\beta^{\sssize(\lower 1pt\hbox{$\ssize 1$}\sssize)}_{k_1\ldots\,k_r}
=\sum_{1\leqslant i<j\leqslant r}
\nabla_{k_i}\!\nabla_{k_j}S\,\cdot\kern -0.8em
\prod\Sb 1\leqslant q\leqslant r\\
\vspace{1pt}
q\neq i,q\neq j\endSb
\kern -0.3em\nabla_{k_q}S\text{\ \ \ for \ \ }r\geqslant 2.
\tag1.16
\endgather
$$
Now let's substitute the expressions \thetag{1.15} and \thetag{1.16}
into the formula \thetag{1.8} and take into account symmetry of
coefficients $a^{k_1\ldots\,k_r}$. This yields
$$
\hskip -2em
\gathered
R_1=\sum^m_{r=1}\sum^n_{k_1=1}\!...\!\sum^n_{k_r\,=1}
r\cdot a^{k_1\ldots\,k_r}\cdot\prod^r_{\rho=2}\nabla_{k_\rho}S
\,\cdot\nabla_{k_1}\,+\\
+\,\sum^m_{r=2}\sum^n_{k_1=1}\!...\!\sum^n_{k_r\,=1}
\frac{r(r-1)}{2}\cdot\nabla_{k_1}\!\nabla_{k_2}S\,\cdot
\prod^r_{\rho=3}\nabla_{k_\rho}S.
\endgathered
\tag1.17
$$
Further we assume that the equations \thetag{1.13} and
\thetag{1.14} with operator \thetag{1.17} in left hand side
are solvable in some sense, and we assume that that $\varphi_{
\sssize(\ssize\lower 1pt\hbox{$\ssize 0$}\sssize)}\neq 0$.
Under these assumptions we concentrate our efforts in studying
the equation \thetag{1.11}.
\head
2. Hamilton-Jacobi equation. 
\endhead
    Note that right hand side of the equality \thetag{1.10} is
similar to that of \thetag{1.1}. This is not casual circumstance.
In order to reveal the nature of this analogy let's denote
$\nabla S$ by $\bold p$. Then $\bold p$ is a covector field with
components
$$
\hskip -2em
p_k=\nabla_kS\text{, \ where \ }k=1,\,\ldots,\,n.
\tag2.1
$$
Using notations \thetag{2.1}, we can rewrite \thetag{1.10} as
$R_0=H(p,\bold p)$, where
$$
\hskip -2em
H(p,\bold p)=\sum^m_{r=0}\sum^n_{k_1=1}\!...\!\sum^n_{k_r\,=1}
a^{k_1\ldots\,k_r}(x^1,\ldots,x^n)\cdot p_{k_1}\cdot\ldots\cdot
p_{k_r}.
\tag2.2
$$
Right hand side of \thetag{2.2} is a polynomial in components
of covector $\bold p$. If $p_1,\,\ldots,\,p_n$ are treated as
independent variables representing components of some arbitrary
covector $\bold p$ at the point $p\in M$, then polynomial $H(p,
\bold p)$ is called the {\bf symbol} of differential operator
\thetag{1.2}. Once $H(p,\bold p)$ is given, the operator
\thetag{1.2} itself can be obtained by substituting $p_k=-i\,
\nabla_k$ into its symbol. This trick is well-known in quantum
mechanics, where momentum vector $\bold p$ is replaced by
differential operator. The analogy with quantum mechanics becomes
more transparent (see \cite{2}) if we take
$$
\hskip -2em
p_k=-\frac{i}{\lambda}\,\nabla_k
\tag2.3
$$
with $\lambda=1/\hbar$ treated as large parameter. If we substitute
\thetag{2.3} into the polynomial \thetag{2.2}, then we obtain
the operator \thetag{1.3} used in the equation \thetag{1.4}.
\par
   Now let's return back to the equation \thetag{1.11}. It is first
equation in the series \thetag{1.9}. Due to \thetag{1.10} and
\thetag{2.2} we can write it as follows:
$$
\hskip -2em
H(p,\nabla S)=0.
\tag2.4
$$
The equation \thetag{2.4} is first order nonlinear partial differential
equation with respect to scalar field $S$. It is often called
{\bf Hamilton-Jacobi equation}. Traditional way of solving such
equations is based on the concept of {\bf characteristic lines}. For
the equation \thetag{2.4} these are parametric curves determined
by {\bf Hamilton equations}:
$$
\xalignat 2
&\hskip -2em
\dot x^i=\frac{\partial H}{\partial p_i},
&&\dot p_i=-\frac{\partial H}{\partial x^i}.
\tag2.5
\endxalignat
$$
In paper \cite{3} the following form of Hamilton equations \thetag{2.5}
was suggested:
$$
\xalignat 2
&\hskip -2em
\dot x^i=\tilde\nabla^i H,
&&\nabla_tp_i=-\nabla_i H.
\tag2.6
\endxalignat
$$
It was especially designed for the case of Hamiltonian dynamical systems
in Riemannian manifolds. Here $\nabla_t$ is standard notation for covariant
differentiation of vector-valued and tensor-valued functions with respect
to parameter $t$ along the curve. In our particular case covariant
derivative $\nabla_tp_i$ is written as
$$
\hskip -2em
\nabla_tp_i=\dot p_i-\sum^n_{k=1}\sum^n_{j=1}\Gamma^k_{ij}
\,p_k\,\dot x^j.
\tag2.7
$$
Covariant derivatives $\nabla_i$ and $\tilde\nabla^i$ in \thetag{2.6}
are less standard objects. The matter is that, writing Hamilton equations
\thetag{2.6}, we treat $H$ as {\bf extended scalar field} in the sense of
the following definition.
\definition{Definition 2.1} {\it Extended} tensor field $\bold X$
of type $(r,s)$ in $M$ is a tensor-valued function that maps each
point $q\in G\subseteq T^*\!M$ to a tensor of the space $T^r_s(p,M)$,
where $p=\pi(q)$. Subset $G$ of $T^*\!M$ is a domain of extended
tensor field $\bold X$. If $G=T^*\!M$, then $\bold X$ is called
{\it global} extended tensor field in $M$. 
\enddefinition
\noindent In local coordinates extended tensor fields are represented
by their components:
$$
\hskip -2em
X^{i_1\ldots\,i_r}_{j_1\ldots\,j_s}=X^{i_1\ldots\,i_r}_{j_1\ldots\,j_s}
(x^1,\ldots,x^n,p_1,\ldots,p_n).
\tag2.8
$$
In contrast to traditional tensor fields, components of extended tensor
field are functions of double set of arguments. First $n$ arguments of
$X^{i_1\ldots\,i_r}_{j_1\ldots\,j_s}$ in \thetag{2.8} are local
coordinates of a point $p\in M$. Others are components of covector
$\bold p\in T^*_p(M)$. Both $p$ and $\bold p$ form a point $q=(p,\bold p)$
of cotangent bundle $T^*\!M$. Further we shall not come deep into the
theory of extended tensor fields, referring reader to paper \cite{3} and
thesis \cite{4} for more details. But we give explicit formulas for
covariant differentiations $\nabla$ and $\tilde\nabla$ in local coordinates
$x^1,\,\ldots,\,x^n,\,p_1,\,\ldots,\,p_n$:
$$
\align
&\hskip -2em
\tilde\nabla^qX^{i_1\ldots\,i_r}_{j_1\ldots\,j_s}=
\frac{\partial X^{i_1\ldots\,i_r}_{j_1\ldots\,j_s}}
{\partial p_q}.
\tag2.9\\
\vspace{2ex}
&\hskip -2em
\aligned
&\nabla_{\!q}X^{i_1\ldots\,i_r}_{j_1\ldots\,j_s}=\frac{\partial
X^{i_1\ldots\,i_r}_{j_1\ldots\,j_s}}{\partial x^q}
+\sum^n_{a=1}\sum^n_{b=1}p_a\,\Gamma^a_{qb}\,\frac{\partial
X^{i_1\ldots\,i_r}_{j_1\ldots\,j_s}}{\partial p_b}\,+\\
&+\sum^r_{k=1}\sum^n_{a_k=1}\!\Gamma^{i_k}_{q\,a_k}\,X^{i_1\ldots\,
a_k\ldots\,i_r}_{j_1\ldots\,\ldots\,\ldots\,j_s}
-\sum^s_{k=1}\sum^n_{b_k=1}\!\Gamma^{b_k}_{q\,j_k}\,
X^{i_1\ldots\,\ldots\,\ldots\,i_r}_{j_1\ldots\,b_k\ldots\,j_s}.
\endaligned
\tag2.10
\endalign
$$
Covariant differentiation $\tilde\nabla$ defined by formula
\thetag{2.9} is called {\bf momentum gradient}. Second covariant
differentiation $\nabla$ defined by formula \thetag{2.10} is
called {\bf spatial gradient}. Note that for traditional tensor
fields formula \thetag{2.10} reduces to standard formula for
covariant derivatives. Therefore use of symbol $\nabla$ in
section~1 and in formula \thetag{2.4} does not contradict to
formula \thetag{2.10}.\par
    Let's take some smooth orientable hypersurface $\sigma$ in $M$.
By $\bold n=\bold n(p)$ we denote smooth field of unit normal vectors
on $\sigma$. Using Riemannian metric $\bold g$ in the manifold $M$,
we can transform it to unit covector field:
$$
\hskip -2em
n_i=\sum^n_{j=1}g_{ij}\,n^j.
\tag2.11
$$
We denote it by the same symbol $\bold n$. Using covector $\bold n=
\bold n(p)$, we set up the following boundary-value problem for
Hamilton-Jacobi equation \thetag{2.4}:
$$
\xalignat 2
&\hskip -2em
S\,\hbox{\vrule height 8pt depth 8pt width 0.5pt}_{\,p\in\sigma}
=s_0=\const,
&&\nabla S\,\hbox{\vrule height 8pt depth 8pt width 0.5pt}_{\,p\in
\sigma}=\nu(p)\cdot\bold n(p).
\tag2.12
\endxalignat
$$
Here $\nu=\nu(p)$ is a scalar factor. It is not arbitrary.
In order to make \thetag{2.12} compatible with the equation
\thetag{2.4} it should satisfy the following condition:
$$
\hskip -2em
H(p,\nu(p)\cdot\bold n(p))=0.
\tag2.13
$$
Aside from boundary-value problem \thetag{2.12} for the equation
\thetag{2.4}, we set up the following initial-value problem for
Hamilton equations \thetag{2.6}:
$$
\xalignat 2
&\hskip -2em
x^i\,\hbox{\vrule height 8pt depth 8pt width 0.5pt}_{\,t=0}
=x^i(p),
&&p_i\,\hbox{\vrule height 8pt depth 8pt width 0.5pt}_{\,t=0}=
\nu(p)\cdot n_i(p).
\tag2.14
\endxalignat
$$
Initial data of Cauchy problem \thetag{2.14} are parametrized by
points of hypersurface $\sigma$. Function $\nu=\nu(p)$ in \thetag{2.14}
is the same scalar factor as in \thetag{2.12}. It is restricted by
condition \thetag{2.13} above.\par
     Solving Cauchy problem \thetag{1.18} for the equations \thetag{2.6},
we get a family of parametric curves in $T^*\!M$. Their projections to
$M$ form a family of parametric curves
$$
\hskip -2em
\gamma=\gamma(t,p)
\tag2.15
$$
starting at the points of hypersurface $\sigma$. Boundary data
\thetag{2.12} and initial data \thetag{2.14} are called {\bf regular}
if function $\nu=\nu(p)$ does not vanish on $\sigma$, i\.\,e\. if
$$
\hskip -2em
\nu(p)\neq 0,
\tag2.16
$$
and if integral curves \thetag{2.15} are transversal to hypersurface
$\sigma$. Now, using Hamilton function $H$, let's define the following
extended scalar field in $M$:
$$
\hskip -2em
\Omega=\sum^n_{i=1}p_i\,\tilde\nabla^iH.
\tag2.17
$$
In terms of $\Omega$ transversality condition $\gamma\nparallel\sigma$
is expressed by inequality
$$
\pagebreak
\hskip -2em
\Omega(p,\nu(p)\cdot\bold n(p))\neq 0.
\tag2.18
$$
Suppose that both regularity conditions \thetag{2.16} and \thetag{2.18}
on hypersurface $\sigma$ are fulfilled. In this case integral curves
\thetag{2.15} fill some neighborhood of hypersurface $\sigma$ and we
can formulate the following theorem.
\proclaim{Theorem 2.1} If boundary data \thetag{2.12} satisfying
condition \thetag{2.13} are regular, i\.\,e\. if both inequalities
\thetag{2.16} and \thetag{2.18} on $\sigma$ are fulfilled, then
boundary-value problem \thetag{2.12} for Hamilton-Jacobi equation
\thetag{2.4} has unique solution $S$ in some neighborhood of
hypersurface $\sigma$. This solution is given by integral
$$
\hskip -2em
S=s_0+\int\limits_\gamma\sum^n_{i=1}p_i\,dx^i
=s_0+\int\limits^{\,t}_{0}\sum^n_{i=1}p_i\,\dot x^i\,dt.
\tag2.19
$$
\endproclaim
\noindent Theorem~2.1 is standard result in the theory of first order
partial differential equations. Its proof can be found in paper \cite{1}.
\head
3. Wave front dynamics. 
\endhead
\parshape 27 0pt 360pt 0pt 360pt 0pt 360pt 0pt 360pt 0pt 360pt
0pt 360pt 0pt 360pt 0pt 360pt 0pt 360pt 0pt 360pt 0pt 360pt
180pt 180pt 180pt 180pt 180pt 180pt 180pt 180pt 180pt 180pt
180pt 180pt 180pt 180pt 180pt 180pt 180pt 180pt 180pt 180pt
180pt 180pt 180pt 180pt 180pt 180pt 180pt 180pt 180pt 180pt
0pt 360pt 
    Theorem~2.1 gives a method for solving Hamilton-Jacobi equation
and finding scalar field $S=S(p)$. Remember that $S$ in \thetag{1.5}
is a phase of propagating wave. By definition {\bf wave front} is a
set of points with constant phase:
$$
\hskip -2em
\sigma(s)=\{p\in M\!:\ S(p)=s=\const\}.
\tag3.1
$$
Hypersurface $\sigma$ in \thetag{2.12} is interpreted as initial
wave front, since $\sigma=\sigma(s_0)$. In regular case for $s$
sufficiently close to $s_0$ wave fronts $\sigma(s)$ in \thetag{3.1}
are smooth hypersurfaces filling some neighborhood of initial
hypersurface $\sigma$. \vadjust{\vskip 32pt\hbox to 0pt{\kern 0pt
\hbox{\special{em:graph pst10e01.gif}}\hss}\vskip -32pt}Curves
\thetag{2.15} are interpreted as rays
in the limit of geometric optics ($\lambda\to\infty$), while Hamilton
equations \thetag{2.6} determine the motion of points, starting at
the time instant $t=0$ and moving along these rays. If we fix their
positions at some nonzero instant of time $t$, for sufficiently small
value of $t$ we would obtain some smooth hypersurface $\sigma_t$
close to initial hypersurface $\sigma$. Hypersurfaces $\sigma_t$
form another family filling some neighborhood of initial wave front
$\sigma$. Should they coincide with $\sigma(s)$\,? In general, the
answer is negative:
$$
\sigma_t\neq\sigma(s),
\tag3.2
$$
see Fig\.~3.1. Hypersurfaces $\sigma_t$ are drawn in solid lines,
while wave fronts $\sigma(s)$ are drawn in dashed lines. Inequality
\thetag{3.2} means that Hamilton equations \thetag{2.6} do not
describe real dynamics of wave fronts exactly. In order to describe
this dynamics in terms of ODE's we should change parameter $t$
in \thetag{2.6} for another parameter $s=S(t)$. Formula \thetag{2.19}
for scalar field $S$ yields the relation of these two parameters:
$$
\hskip -2em
\frac{ds}{dt}=\sum^n_{i=1}p_i\,\dot x^i=\sum^n_{i=1}p_i\,
\tilde\nabla^iH=\Omega.
\tag3.3
$$
Changing $t$ for $s$ according to the formula \thetag{3.3}, we
keep symbol $t$ for denoting this new parameter. Then from
\thetag{2.6} we derive modified Hamilton equations
$$
\xalignat 2
&\hskip -2em
\dot x^i=\frac{\tilde\nabla^i H}{\Omega},
&&\nabla_tp_i=-\frac{\nabla_i H}{\Omega},
\tag3.4
\endxalignat
$$
where $\Omega$ is given by formula \thetag{2.17}. Differential
equations \thetag{3.4} are called {\bf differential equations
of wave front dynamics}.
\head
4. Legendre transformation.
\endhead
    Hamilton equations \thetag{2.6} and modified Hamilton equations
\thetag{3.4} are the equations in cotangent bundle $T^*\!M$. Now we
shall transform them into the equations in tangent bundle $TM$. Note
that for each point $q=(p,\bold p)$ of cotangent bundle $T^*\!M$
we can consider vector $\bold v\in T_p(M)$ with components
$$
\hskip -2em
v^i=\tilde\nabla^i H.
\tag4.1
$$
Joining together $p$ and $\bold v$, we get a point $q=(p,\bold v)$ of
tangent bundle $TM$. This determines a map $\lambda^{-1}\!: T^*\!M\to
TM$, which is known as {\bf inverse Legendre transformation}. We shall
assume it to be invertible and denote by $\lambda$ direct Legendre
transformation: $\lambda\!: TM\to T^*\!M$. Suppose that $\bold X$ is
some extended tensor field. According to the definition~2.1, it is
tensor-valued function in $T^*\!M$. Then composition $\bold X\compos
\lambda$ is tensor-valued function in tangent bundle $TM$. So we
get another definition of extended tensor field.
\definition{Definition 4.1} {\it Extended} tensor field $\bold X$
of type $(r,s)$ in $M$ is a tensor-valued function that maps each
point $q\in G\subseteq TM$ to a tensor of the space $T^r_s(p,M)$,
where $p=\pi(q)$. Subset $G$ of $TM$ is a domain of extended
tensor field $\bold X$. If $G=TM$, then $\bold X$ is called
{\it global} extended tensor field in $M$. 
\enddefinition
    In the case of arbitrary smooth manifold $M$ definitions~2.1 and
4.1 lead to different theories. But for Riemannian manifold $M$ tangent
bundle $TM$ and cotangent bundle $T^*\!M$ are bound with each other
by duality maps
$$
\xalignat 2
&\hskip -2em
\bold g\!:TM\to T^*\!M,
&&\bold g^{-1}\!:T^*\!M\to TM.
\tag4.2
\endxalignat
$$
In local coordinates duality maps \thetag{4.2} are represented as
index lowering and index raising procedures in arguments of $\bold X$.
Due to duality maps \thetag{4.2} two objects introduced by
definitions~2.1 and 4.1 are the same in essential. We call them
{\it covariant} and {\it contravariant} representations of extended
tensor field $\bold X$. If we take contravariant representation of
$\bold X$, then formulas \thetag{2.9} and \thetag{2.10} are rewritten as
$$
\align
&\hskip -2em
\tilde\nabla_qX^{i_1\ldots\,i_r}_{j_1\ldots\,j_s}=\frac{\partial
X^{i_1\ldots\,i_r}_{j_1\ldots\,j_s}}{\partial p^q},
\tag4.3\\
\vspace{2ex}
&\hskip -2em
\aligned
&\nabla_qX^{i_1\ldots\,i_r}_{j_1\ldots\,j_s}=\frac{\partial
X^{i_1\ldots\,i_r}_{j_1\ldots\,j_s}}{\partial x^q}
-\sum^n_{a=1}\sum^n_{b=1}p^a\,\Gamma^b_{qa}\,\frac{\partial
X^{i_1\ldots\,i_r}_{j_1\ldots\,j_s}}{\partial p^b}\,+\\
&+\sum^r_{k=1}\sum^n_{a_k=1}\!\Gamma^{i_k}_{q\,a_k}\,X^{i_1\ldots\,
a_k\ldots\,i_r}_{j_1\ldots\,\ldots\,\ldots\,j_s}
-\sum^s_{k=1}\sum^n_{b_k=1}\!\Gamma^{b_k}_{q\,j_k}
X^{i_1\ldots\,\ldots\,\ldots\,i_r}_{j_1\ldots\,b_k\ldots\,j_s}.
\endaligned
\tag4.4
\endalign
$$
Direct and inverse Legendre transformations are other two maps
$$
\xalignat 2
&\hskip -2em
\lambda\!:TM\to T^*\!M,
&&\lambda^{-1}\!:T^*\!M\to TM
\tag4.5
\endxalignat
$$
binding $TM$ and $T^*\!M$. Presence of nonlinear maps \thetag{4.5}
increases the number of representations of extended tensor field
$\bold X$. If $\bold Y=\bold X\compos\lambda^{-1}$, then we
say that $\bold Y$ is $\bold p$-representation or {\bf momentum
representation} for $\bold X$, while $\bold X$ is called\linebreak
$\bold v$-representation or {\bf velocity representation} for $\bold Y$.
Total set of representations for extended tensor field is given
in the following table:
$$
\vbox{\offinterlineskip
\def\vr{\vrule height 14pt depth 8pt}
\settabs\+\hskip 0.5cm&\hskip 4.5cm&\hskip 0.5cm&\hskip 5.2cm&\cr
\hrule
\+\vr\strut &covariant $\bold p$-representation
 &\vr\strut &contravariant $\bold p$-representation &\vr\strut\cr
\hrule
\+\vr\strut &covariant $\bold v$-representation
 &\vr\strut &contravariant $\bold v$-representation &\vr\strut\cr
\hrule}
$$
The following extended scalar field is called {\bf Lagrange function}:
$$
\hskip -2em
l=\sum^n_{i=1}p_i\cdot\tilde\nabla^iH-H
\tag4.6
$$
Covariant $\bold p$ representation is natural for Hamilton
function $\bold H$, while Lagrange function is naturally
used in contravariant $\bold v$-representation: $L=l\compos\lambda$.
If Lagrange function \thetag{4.6} is already transformed to
contravariant $\bold v$-representation, then direct Legendre
transformation $\lambda$ in \thetag{4.5} can be given by formula
similar to \thetag{4.1}:
$$
\hskip -2em
p_i=\tilde\nabla_{\!i}L
\tag4.7
$$
(see \cite{3} for more details). Note that formulas \thetag{2.9},
\thetag{2.10}, \thetag{4.3}, and \thetag{4.4} for $\nabla$ and
$\tilde\nabla$ are valid either in $\bold p$-representation and
in $\bold v$-representation. However, in \linebreak
$\bold v$-representation
we should replace $p_i$ by $v_i$ and $p^i$ by $v^i$ in these
formulas. Therefore $\tilde\nabla$ is called {\bf velocity gradient}
in $\bold v$-representation.\par
    Duality maps \thetag{4.2} commutate with covariant differentiations
$\nabla$ and $\tilde\nabla$. This is expressed by the following
equalities:
$$
\xalignat 2
&\hskip -2em
\nabla(\bold X\compos\bold g)=(\nabla\bold X)\compos\bold g,
&&\tilde\nabla(\bold X\compos\bold g)=(\tilde\nabla\bold X)
\compos\bold g.
\tag4.8
\endxalignat
$$
Legendre transformation $\lambda$ does not commutate with differentiations
$\nabla$ and $\tilde\nabla$. Unlike \thetag{4.8}, we have the
following equalities, where $C$ is the operation of contraction:
$$
\align
&\hskip -2em
\tilde\nabla(\bold Y\compos\lambda)=C(\tilde\nabla\bold Y\compos
\lambda\otimes\tilde\nabla\tilde\nabla L),
\tag4.9\\
&\hskip -2em
\nabla(\bold Y\compos\lambda)=\nabla\bold Y\compos\lambda
+C(\tilde\nabla\bold Y\compos\lambda\otimes\nabla\tilde\nabla L).
\tag4.10
\endalign
$$
In local coordinates the equalities \thetag{4.9} and \thetag{4.10}
are written as follows:
$$
\align
&\hskip -2em
\tilde\nabla_{\!q}X^{i_1\ldots\,i_r}_{j_1\ldots\,j_s}
=\sum^n_{k=1}\tilde\nabla_{\!q}\tilde\nabla_kL
\cdot\tilde\nabla^kY^{i_1\ldots\,i_r}_{j_1\ldots\,j_s},
\tag4.11\\
&\hskip -2em
\nabla_{\!q}X^{i_1\ldots\,i_r}_{j_1\ldots\,j_s}
=\nabla_{\!q}Y^{i_1\ldots\,i_r}_{j_1\ldots\,j_s}+
\sum^n_{k=1}\nabla_{\!q}\tilde\nabla_kL
\cdot\tilde\nabla^kY^{i_1\ldots\,i_r}_{j_1\ldots\,j_s}.
\tag4.12
\endalign
$$
Similarly, one can write inverse relationships for \thetag{4.9} and
\thetag{4.10}:
$$
\align
&\hskip -2em
\tilde\nabla(\bold X\compos\lambda^{-1})=C(\tilde\nabla\bold X\compos
\lambda^{-1}\otimes\tilde\nabla\tilde\nabla H),
\tag4.13\\
&\hskip -2em
\nabla(\bold X\compos\lambda^{-1})=\nabla\bold X\compos\lambda^{-1}
+C(\tilde\nabla\bold X\compos\lambda^{-1}\otimes\nabla\tilde\nabla H).
\tag4.14
\endalign
$$
In local coordinates these relationships \thetag{4.13} and \thetag{4.14}
are expressed as
$$
\align
&\hskip -2em
\tilde\nabla^kY^{i_1\ldots\,i_r}_{j_1\ldots\,j_s}
=\sum^n_{q=1}\tilde\nabla^k\tilde\nabla^qH
\cdot\tilde\nabla_{\!q}X^{i_1\ldots\,i_r}_{j_1\ldots\,j_s},
\tag4.15\\
&\hskip -2em
\nabla_kY^{i_1\ldots\,i_r}_{j_1\ldots\,j_s}
=\nabla_kX^{i_1\ldots\,i_r}_{j_1\ldots\,j_s}+
\sum^n_{q=1}\nabla_k\tilde\nabla^qH
\cdot\tilde\nabla_{\!q}X^{i_1\ldots\,i_r}_{j_1\ldots\,j_s}.
\tag4.16
\endalign
$$
In \thetag{4.11}, \thetag{4.12}, \thetag{4.15}, and \thetag{4.16}
we assume that $\bold X=\bold Y\compos\lambda$ and $\bold Y=\bold X
\compos\lambda^{-1}$.\par
    Now let's consider Lagrange function $L$. Its
$\bold p$-representation is given by formula \thetag{4.6}. Let's
calculate $\nabla_kl$ directly by means of formula \thetag{4.6}:
$$
\hskip -2em
\nabla_kl=\sum^n_{i=1}p_i\cdot\nabla_k\tilde\nabla^iH-\nabla_kH.
\tag4.17
$$
Then let's apply formula \thetag{4.16} and let's calculate this
derivative $\nabla_kl$ again:
$$
\hskip -2em
\nabla_kl=\nabla_kL+\sum^n_{q=1}\nabla_k\tilde\nabla^qH\cdot
\tilde\nabla_qL.
\tag4.18
$$
Comparing \thetag{4.17} and \thetag{4.18}, and taking into account
\thetag{4.7}, we get:
$$
\hskip -2em
\nabla_kL=-\nabla_kH.
\tag4.19
$$
In coordinate free form this equality \thetag{4.19} looks like
$\nabla L=-\nabla H\compos\lambda$.
\head
5. Lagrangian representation for
the equations of wave front dynamics.
\endhead
    Now we are ready to transform differential equations \thetag{3.4}
to $\bold v$-representation. Using formula \thetag{4.1}, we write fist
equation \thetag{3.4} as
$$
\hskip -2em
\dot x^i=\frac{v^i}{\Omega}.
\tag5.1
$$
Extended scalar field $\Omega$ in $\bold p$-representation is
determined by formula \thetag{2.17}. Due to \thetag{4.1} and
\thetag{4.7} its $\bold v$-representation is determined by formula
$$
\hskip -2em
\Omega=\sum^n_{i=1}v^i\cdot\tilde\nabla_i L.
\tag5.2
$$
In order to transform second equation \thetag{3.4} we use formulas
\thetag{4.7} and \thetag{4.19}:
$$
\pagebreak
\hskip -2em
\nabla_t(\tilde\nabla_iL)=\frac{\nabla_iL}{\Omega}.
\tag5.3
$$
Written together, \thetag{5.1} and \thetag{5.3} form Lagrangian
representation for the equations of wave front dynamics \thetag{3.4}.
Here are these equations:
$$
\xalignat 2
&\hskip -2em
\dot x^i=\frac{v^i}{\Omega},
&&\nabla_t(\tilde\nabla_iL)=\frac{\nabla_iL}{\Omega}.
\tag5.4
\endxalignat
$$
Though written in terms of Lagrange function, the equations
\thetag{5.4} are not Euler-Lagrange equations due to the
denominator $\Omega$ in them.
\head
6. Dynamical systems admitting normal shift.
\endhead
    Note that vector $\bold v$ with components $v^1,\,\ldots,\,v^n$
is not an actual velocity vector for the dynamics described by
differential equations \thetag{5.4}. Let's denote actual velocity
vector by $\bold u$. Then from first equation \thetag{5.4} we derive
$$
\hskip -2em
u^i=\frac{v^i}{\Omega}.
\tag6.1
$$
Denominator $\Omega$ in \thetag{6.1} is not constant, it depends
on $x^1,\,\ldots,\,x^n$ and on components of velocity vector $\bold v$
as well. Therefore \thetag{6.1} defines nonlinear map
$$
\hskip -2em
\mu\!: TM\to TM
\tag6.2
$$
similar to Legendre map $\lambda$ determined by \thetag{4.7}. If
this map \thetag{6.2} invertible, then one can transform \thetag{5.4}
to Newtonian form:
$$
\xalignat 2
&\hskip -2em
\dot x^i=u^i,
&&\nabla_tu^i=F^i(x^1,\ldots,x^n,u^1,\ldots,u^n).
\tag6.3
\endxalignat
$$
Components of force vector $\bold F$ in \thetag{6.3} are expressed
through Lagrange function $L$ and its derivatives. Now we shall not
derive general formula for $F^i$. But we consider very important special
case, when $L$ is fiberwise spherically symmetric extended scalar
field in Riemannian manifold $M$.
\definition{Definition 6.1} Extended tensor field $\bold X$ is called
{\it fiberwise spherically symmetric} if its components depend only
on modulus of velocity vector:
$$
X^{i_1\ldots\,i_r}_{j_1\ldots\,j_s}=X^{i_1\ldots\,i_r}_{j_1\ldots\,
j_s}(x^1,\,\ldots,\,x^n,v)\text{, \ where \ }v=|\bold v|.
$$
This means that $X^{i_1\ldots\,i_r}_{j_1\ldots\,j_s}$ is spherically
symmetric function within each fiber of tangent bundle $TM$ for each
fixed point $p\in M$.
\enddefinition
    Fiberwise spherically symmetric tensor fields were considered
in Chapter~\uppercase\expandafter{\romannumeral 7} of thesis \cite{4},
see also recent papers \cite{3} and \cite{5}. Now suppose that $L$ is
fiberwise spherically symmetric extended scalar field:
$L=L(x^1,\ldots,x^n,v)$. Then from \thetag{5.2} we derive that
$\Omega$ is also fiberwise spherically symmetric scalar field:
$$
\hskip -2em
\Omega=|\bold v|\cdot L'.
\tag6.4
$$
Here $L'$ is partial derivative of $L$ with respect of its last
argument $v=|\bold v|$. Substituting \thetag{6.4} into the equality
\thetag{6.1}, we derive the following formulas:
$$
\xalignat 2
&\hskip -2em
u^i=\frac{N^i}{L'},
&&|\bold u|=\frac{1}{|L'|}.
\tag6.5
\endxalignat
$$
By $N^i$ in \thetag{6.5} we denote components of unit vector $\bold N$
directed along vector $\bold v$. Using its components, we define
two operators of orthogonal projection $\bold Q$ and $\bold P$:
$$
\xalignat 2
&\hskip -2em
Q^i_k=\frac{v^i\,v_k}{|\bold v|^2},
&&P^i_k=\delta^i_k-\frac{v^i\,v_k}{|\bold v|^2}.
\tag6.6
\endxalignat
$$
First of them is a projector to the direction of vector $\bold v$,
second is a projector to hyperplane perpendicular to $\bold v$.
Projection operators $\bold Q$ and $\bold P$ with components
\thetag{6.6} are complementary to each other, this means that
$$
\bold Q+\bold P=\bold 1\quad\text{and}\quad\bold Q\compos\bold P
=\bold P\compos\bold Q=\bold 0.
$$
Now let's consider second equation \thetag{5.4}. For $\nabla_t(\tilde
\nabla_iL)$ in it we have
$$
\hskip -2em
\nabla_t(\tilde\nabla_iL)=
\sum^n_{k=1}\nabla_k\tilde\nabla_iL\cdot\dot x^k
+\sum^n_{k=1}\tilde\nabla_k\tilde\nabla_iL\cdot
\nabla_t v^k.
\tag6.7
$$
By direct calculations for various derivatives in \thetag{6.7} we obtain
$$
\align
&\hskip -2em
\tilde\nabla_iL=L'\cdot\frac{v_i}{|\bold v|}=L'\cdot N_i,
\tag6.8
\\
&\hskip -2em
\nabla_k\!\tilde\nabla_iL=\nabla_kL'\cdot N_i,
\tag6.9
\\
&\hskip -2em
\tilde\nabla_k\!\tilde\nabla_iL=L''\cdot Q_{ik}
+\frac{L'}{|\bold v|}\cdot P_{ik}.
\tag6.10
\endalign
$$
Substituting \thetag{6.9} and \thetag{6.10} into \thetag{6.7},
we get the following equality:
$$
\hskip -2em
\nabla_t(\tilde\nabla_iL)=\sum^n_{k=1}\frac{\nabla_kL'}{L'}\cdot
Q^k_i+\sum^n_{k=1}\left(L''\cdot Q_{ik}+\frac{L'}{|\bold v|}\cdot
P_{ik}\right)\cdot\nabla_tv^k.
\tag6.11
$$
Next step is to relate $\nabla_tv^k$ in \thetag{6.11} with
$\nabla_tu^k$. For this purpose let's apply $\nabla_t$ to
the first equality in \thetag{6.5}. As a result of direct
calculations we get:
$$
\nabla_tu^k=\sum^n_{s=1}\left(\frac{1}{|\bold v|\cdot L'}
\cdot P^k_s-\frac{L''}{(L')^2}\cdot Q^k_s\right)\cdot\nabla_tv^s
-\sum^n_{s=1}\frac{\nabla_{\!s}L'}{(L')^3}\cdot Q^{sk}.
\tag6.12
$$
In order to express $\nabla_tv^s$ through $\nabla_tu^k$ we should
invert matrix $D$ with components
$$
D^k_s=\frac{1}{|\bold v|\cdot L'}\cdot P^k_s
-\frac{L''}{(L')^2}\cdot Q^k_s.
$$
This matrix is invertible if and only if $L''\neq 0$. The condition
$L'\neq 0$ is fulfilled since $\Omega=|\bold v|\cdot L'\neq 0$. This
follows from transversality condition \thetag{2.18}. Inverting matrix
$D$ and applying it to \thetag{6.12}, we get the following expression
for $\nabla_tv^k$:
$$
\nabla_tv^k=\sum^n_{s=1}\left(|\bold v|\cdot L'\cdot P^k_s
-\frac{(L')^2}{L''}\cdot Q^k_s\right)\cdot\nabla_tu^s
-\sum^n_{s=1}\frac{\nabla_{\!s}L'}{L'\cdot L''}\cdot Q^{sk}.
\tag6.13
$$
Now let's substitute \thetag{6.13} back to right hand side of
\thetag{6.11}. As a result we get
$$
\hskip -2em
\nabla_t(\tilde\nabla_iL)=(L')^2\sum^n_{k=1}(P_{ik}-Q_{ik})\cdot
\nabla_tu^k.
\tag6.14
$$
Then let's substitute \thetag{6.14} into second equation \thetag{5.4}
and let's use formula \thetag{6.4} for denominator $\Omega$ in it.
This leads to the equation for $\nabla_tu^k$:
$$
\hskip -2em
(L')^2\sum^n_{k=1}(P_{ik}-Q_{ik})\cdot\nabla_tu^k=
\frac{\nabla_iL}{|\bold v|\cdot L'}.
\tag6.15
$$
Equation \thetag{6.15} can be explicitly solved. For $\nabla_tu^k$
from \thetag{6.15} we derive
$$
\hskip -2em
\nabla_tu^k=\sum^n_{i=1}\frac{\nabla_iL}{|\bold v|\cdot (L')^3}\cdot
\bigl(g^{ik}-2\,N^i\,N^k\bigr).
\tag6.16
$$
Comparing \thetag{6.16} with \thetag{6.3} we see that right hand side
of \thetag{6.16} is the expression for $F^k$. Formula for covariant
components $F_k$ is more elegant:
$$
\hskip -2em
F_k=\sum^n_{i=1}\frac{\nabla_iL}{|\bold v|\cdot (L')^3}\cdot
\left(\delta^i_k-2\,N^i\,N_k\right).
\tag6.17
$$
The only problem now is that $\nabla_iL$, $L'$, and $|\bold v|$
in right hand side of \thetag{6.17} are given in $\bold v$-representation,
while components of force vector in \thetag{6.3} should depend on
components of actual velocity $\bold u$. Our last effort is to transform
\thetag{6.17} to\linebreak $\bold u$-representation, using map
\thetag{6.2}.\par
    Let's denote by $h$ Hamilton function $H$ in $\bold v$-representation,
i\.\,e\. $h=H\compos\lambda$. Then for $h$ we have the following formula
similar to formula \thetag{4.6} for $l$:
$$
\hskip -2em
h=\sum^n_{i=1}v^i\cdot\tilde\nabla_iL-L=|\bold v|\cdot L'-L.
\tag6.18
$$
Differentiating \thetag{6.18}, we obtain
$$
\xalignat 2
&\hskip -2em
\nabla_ih=v\cdot\nabla_i L'-\nabla_iL,
&&h'=v\cdot L''.
\tag6.19
\endxalignat
$$
Let $W$ be $\bold u$-representation for $h$. Denote by $u$ modulus of
actual velocity vector. From \thetag{6.5} and \thetag{6.18} we conclude
that $W$ is a fiberwise spherically symmetric scalar field in
$\bold u$-representation: $W=W(x^1,\ldots,x^n,u)$, where $u=|\bold u|$.
Now let's calculate derivatives of $W$. For $W'$ and $\nabla_iW$
we get:
$$
\hskip -2em
\aligned
&W'=\frac{\partial W}{\partial u}=
\frac{\partial h}{\partial v}\cdot\frac{\partial v}{\partial u}=
h'\cdot\fracwithdelims(){\partial u}{\partial v}^{\lower 2pt
\hbox{$\ssize\!-1$}},\\
\vspace{2ex}
&\nabla_iW=\frac{\partial W}{\partial x^i}=
\nabla_ih-h'\cdot\frac{\partial u}{\partial x^i}\cdot
\fracwithdelims(){\partial u}{\partial v}^{\lower 2pt
\hbox{$\ssize\!-1$}}.
\endaligned
\tag6.20
$$
Now let's remember that $u$ and $v$ are related by second equality
in \thetag{6.5}. We write it as $u=\varepsilon/L'$, where $\varepsilon
=\pm 1$ depending on sign of $L'$. Then we get
$$
\xalignat 2
&\hskip -2em
\frac{\partial u}{\partial v}=-\frac{\varepsilon\cdot L''}{(L')^2},
&&\frac{\partial u}{\partial x^i}=-\frac{\varepsilon\cdot\nabla_iL'}
{(L')^2}.
\tag6.21
\endxalignat
$$
Combining \thetag{6.19}, \thetag{6.20}, and \thetag{6.21}, we
can transform \thetag{6.22} as follows:
$$
\xalignat 2
&\hskip -2em
W'=-\varepsilon\cdot|\bold v|\cdot (L')^2,
&&\nabla_iW=-\nabla_iL.
\tag6.22
\endxalignat
$$
If we write second equality \thetag{6.5} as $|\bold u|=\varepsilon/L'$,
then from \thetag{6.22} we derive:
$$
\hskip -2em
|\bold u|\cdot\frac{\nabla_iW}{W'}=\frac{\nabla_iL}{|\bold v|\cdot(L')^3}.
\tag6.23
$$
Comparing the equality \thetag{6.23} with formula \thetag{6.17} for
components of force vector $\bold F$, we can rewrite this formula
in the following form: 
$$
\hskip -2em
F_k=-|\bold u|\cdot\sum^n_{i=1}\frac{\nabla_iW}{W'}\cdot
\left(2\,N^i\,N_k-\delta^i_k\right).
\tag6.24
$$
Now formula \thetag{6.24} is completely compatible with \thetag{6.3}.
Its right hand side is given in $\bold u$-representation as it is
required in \thetag{6.3}.\par
    Formula determines force field of Newtonian dynamical system
\thetag{6.3} describing the dynamics of wave front for linear
wave equation \thetag{1.4}. But most important fact for us is that
formula \thetag{6.24} gives the link to the theory of Newtonian
dynamical systems admitting normal shift of hypersurfaces. This
theory was developed in series of papers \cite{6--21}, on the base
of which theses \cite{4} and \cite{22} were prepared. According to
results of Chapter~\uppercase\expandafter{\romannumeral 7} of thesis
\cite{4} (see also paper \cite{5}), force field of any Newtonian
dynamical system admitting normal shift of hypersurfaces in
Riemannian manifold of the dimension $n\geqslant 3$ is given by
explicit formula:
$$
\hskip -2em
F_k=\frac{h(W)\cdot N_k}{W'}
-|\bold u|\cdot\sum^n_{i=1}\frac{\nabla_iW}{W'}\cdot
\left(2\,N^i\,N_k-\delta^i_k\right).
\tag6.25
$$
Here $W=W(x^1,\ldots,x^n,|\bold u|)$ is arbitrary fiberwise
spherically symmetric extended scalar field with $W'\neq 0$,
and $h=h(w)$ is arbitrary function of one variable.
\head
7. Conclusions.
\endhead
    Comparing formulas \thetag{6.24} and \thetag{6.25}, we see
that they do coincide for the case $h=0$. Term with $h=h(W)$
in \thetag{6.25} is responsible for energy dissipation and
energy pumping phenomena. In the absence of these phenomena,
i\.\,e\. for $h=0$, we can state the following main results
of present paper:
\roster
\rosteritemwd=1pt
\item"1)" conservative dynamical systems admitting normal shift of
hypersurfaces in Riemannian manifold of the dimension $n\geqslant 3$
coincide with those describing wave front dynamics in quasiclassical
limit for wave operators with fiberwise spherically symmetric symbol
$H$;
\item"2)" these systems are not Hamiltonian, but they are very close
to Hamiltonian systems (see equations \thetag{3.4} above);
\item"3)" they possess first integral, which can be interpreted as
energy (or as Hamilton function in $\bold p$-representation).
\endroster
Third result that $W$ is a first integral of Newtonian dynamical system
with force field \thetag{6.24} was already known (see paper \cite{5} or
Chapter~\uppercase\expandafter{\romannumeral 7} of thesis \cite{4}). But
its interpretation was not so clear as it is now.\par
    Note that the equations of wave front dynamics do not exhaust the
whole class of Newtonian dynamical systems admitting normal shift of
hypersurfaces in $M$. The problem of proper interpretation for
non-conservative term
$$
\frac{h(W)\cdot N_k}{W'}
$$
in formula \thetag{6.25} in the sense of theory from \cite{1} and 
\thetag{23} is still open. This problem should be studied in separate
paper. Methods developed in book \cite{24} might appear useful for
this purpose.\par
   {\bf Remark}. Technique based on usage of extended tensor fields is
not new (see for instance book \cite{25}). And, as noted by Prof\.~Mircea
Crasmareanu, our terminology above is not standard. However, we do
not change this terminology in order to keep integrity of present
paper and our previous papers (see \cite{3-22} and some other
papers in Electronic Archive at LANL\footnotemark).
\footnotetext{Electronic Archive at Los Alamos National Laboratory of USA
(LANL). Archive is accessible through Internet 
{\bf http:/\negskp/arXiv.org}, it has mirror site 
{\bf http:/\negskp/ru.arXiv.org} at the Institute for Theoretical and
Experimental Physics (ITEP, Moscow).}

\enddocument
\end